# SPECTRAL GAP FOR THE ZERO RANGE PROCESS WITH CONSTANT RATE

By Ben Morris

*University of California, Davis*

We solve an open problem concerning the relaxation time (inverse spectral gap) of the zero range process in $\mathbf{Z}^d/L\mathbf{Z}^d$ with constant rate, proving a tight upper bound of $O((\rho+1)^2 L^2)$, where $\rho$ is the density of particles.

## 1. Introduction.

1.1. *Basic definitions.* Letting $G = (V, E)$ be a finite, connected, regular graph, we define a *configuration* as follows. In a configuration, a total of $r$ indistinguishable particles are distributed over the vertices in $V$. For a configuration $\eta$, we denote by $\eta(v)$ the number of particles in vertex $v$ [so that $\sum_{v \in V} \eta(v) = r$]. We define the density of particles as $\rho = r/|V|$.

The *zero range process with rate* $\lambda$ is the following continuous-time Markov process on configurations. Suppose that the current state is $\eta$. For each vertex $v$ at rate $\lambda$, if $\eta(v) > 0$, we expel a particle from $v$ to a random neighbor, that is, decrease $\eta(v)$ by one, choose a neighbor $w$ of $v$ uniformly at random (u.a.r.) and increase $\eta(w)$ by one.

Note that since the process is irreducible and has symmetric transition rates, the distribution at time $t$ converges to uniform as $t \to \infty$. Let $\mathcal{C}$ denote the space of configurations, and for probability distributions $\mu, \nu$ on $\mathcal{C}$, let

$$\|\mu - \nu\| = \max_{Q \subset \mathcal{C}} \mu(Q) - \nu(Q) = \min_{X \sim \mu, Y \sim \nu} \mathbf{P}(X \neq Y)$$

be the total variation distance. The *spectral gap*, defined as the absolute value of the second largest eigenvalue of the generator of the process, governs the asymptotic rate of convergence to the stationary distribution (see, e.g., [6]). More specifically, we have

(1) $$\text{gap} = \min_x \lim_{t \to \infty} -\frac{1}{t} \log \|K^t(x, \cdot) - \mathcal{U}\|,$$

where $K$ is the transition kernel for the zero range process and $\mathcal{U}$ is the uniform distribution over configurations. We define the *relaxation time* to be

$$\tau_{\text{relax}} = 1/\text{gap} = \max_x \lim_{t \to \infty} -t/\log \|K^t(x, \cdot) - \mathcal{U}\|.$$

Thus, the relaxation time is the smallest value of $\tau$ such that for every starting configuration $x$, there is a constant $C$ such that $\|K^t(x, \cdot) - \mathcal{U}\| \leq Ce^{-t/\tau}$ for all $t \geq 0$.

This paper is concerned with bounding the relaxation time in the important special case where $G$ is the $d$-dimensional torus $\mathbf{Z}^d/L\mathbf{Z}^d$. Note that the spectral gap is proportional to the rate $\lambda$, so we can assume without loss of generality that $\lambda = 1$. We will call this process the *ZRP on G*. We will take the dimension $d$ to be arbitrary but fixed and bound the relaxation time as a function of $\rho$ and $L$.

1.2. *Background, motivation and summary of results.* The zero range process is a widely studied Markov chain in statistical mechanics. In the general zero range model, there is a rate function $c(k)$ that gives the rate at which a site occupied by $k$ particles expels one [e.g., the assumption of independent random walks corresponds to $c(k) = k$]. In our model, we have $c(k) = 1$ for all $k$, so the rate at which a vertex expels a particle does not depend on the number of particles there. In the language of queuing theory, this is a closed queuing network with i.i.d. exponential service times.

A key estimate needed for hydrodynamic limits is a sharp upper bound on the relaxation time (see [11]). Landim, Sethuraman and Varadhan [5] used the techniques of Lu and Yau [7] to prove that under certain conditions on the rate function $c(\cdot)$, the relaxation time is $O(L^2)$; see also [1]. However, their conditions did not apply to a bounded rate function since they required that there be a constant $\gamma > 0$ such that $c(k) > \gamma k$ for sufficiently large $k$. Indeed, as they observe (see Example 1.1 of their paper), their result does not apply to the case of a constant rate since their bound holds uniformly in $\rho$, but the relaxation time is of order $(1+\rho)^2 L^2$ when the rate is constant and $d=1$. In fact, a lower bound of order $(1+\rho)^2 L^2$ holds for all $d$; this can be seen by substituting the test function $f(\eta) = \sum_{v \in V} \eta(v) \cos(\frac{\pi v_1}{L})$ into the variational definition of spectral gap (see, e.g., [4]). This type of function was used by David Wilson in [12].

While we cannot expect an $O(L^2)$ bound to hold uniformly in $\rho$, it is natural to ask whether such a bound holds if we incorporate the extra factor $(\rho + 1)^2$. In this paper, we prove that this is indeed the case.

We will study the ZRP on the torus $\mathbf{Z}^d/L\mathbf{Z}^d$, but our proof is easily extended to variants such as the box $\{0, 1, \ldots, L\}^d$. Indeed, our proof shows that for an arbitrary regular graph $G$, the bound on the spectral gap for random walk on $G$ given by the paths technique (see [3, 9]) also applies, up to constant factors, to the ZRP on $G$. Similar results were already known for the exclusion process; see [2].

**2. Comparison technique.** It turns out that the ZRP is easier to analyze on the complete graph than on the torus. Fortunately, comparison techniques (see, e.g., [2]) allow us to reduce to this case. Let $G_1 = \mathbf{Z}^d/L\mathbf{Z}^d$ and let $G_2$ be the complete graph on $L^d$ vertices (not including self-loops). Let $\mathcal{V}_1$ denote the set of configurations of the ZRP on $G_1$, let $\mathcal{E}_1 = \{(\eta, \eta'): \eta, \eta' \in \mathcal{V}_1$ and $\eta \to \eta'$ is a transition in the ZRP on $G_1\}$ and let $\mathcal{G}_1$ be the graph $\langle \mathcal{V}_1, \mathcal{E}_1 \rangle$. Note that the ZRP on $G_1$ can be described as the continuous-time random walk on $\mathcal{G}_1$ in which edges are crossed at rate $1/2d$. Similarly, if we define $\mathcal{V}_2, \mathcal{E}_2$ and $\mathcal{G}_2$ analogously (with $G_2$ replacing $G_1$), then the ZRP on $G_2$ is the random walk on $\mathcal{G}_2$ in which edges are crossed at rate $1/(L^d-1)$. Fairly standard techniques allow us to estimate the relaxation time for the ZRP on $G_1$ by comparing it with the ZRP on $G_2$. More specifically, we will set up a suitable multicommodity flow on $\mathcal{G}_1$ in which we route one unit of flow from $\eta$ to $\eta'$ for each $(\eta, \eta') \in \mathcal{E}_2$ simultaneously. For any such flow $f$ and any edge $e$ in $\mathcal{E}_1$, let $f(e)$ denote the total flow along $e$; that is, $f(e)$ is the sum over all ordered pairs $(\eta, \eta') \in \mathcal{E}_2$ of the $\eta \to \eta'$ flow carried by $e$. Define the *congestion* $\mathcal{C}(f) = \frac{1}{L^d-1} \max_e f(e)$, that is, the maximum flow along any edge normalized by $L^d - 1$, and the *length* $\mathcal{L}(f)$ to be the length of a longest flow-carrying path. The following theorem is a special case of results in [2] (see also [10]):

THEOREM 1 ([2]). *Let $\tau_1$ and $\tau_2$ denote the relaxation times for the ZRP's on $G_1$ and $G_2$, respectively. For any flow $f$, we have $\tau_1 \leq 2d\mathcal{C}(f)\mathcal{L}(f)\tau_2$.*

We will bound $\tau_1$ by constructing a flow $f$ with congestion $\mathcal{C}(f) \leq L$ and length $\mathcal{L}(f) \leq dL$. By Theorem 1, this implies $\tau_1 \leq 2d^2 L^2 \tau_2$. We will now describe the flow $f$. Let $(\eta, \eta')$ be an edge in $\mathcal{E}_2$. Then the configurations $\eta$ and $\eta'$ differ at only two vertices $u$ and $v$, with, say, $\eta(u) = \eta'(u) + 1$ and $\eta(v) = \eta'(v) - 1$. For a vertex $w$, define $\chi_w$ to be the configuration that has just a single particle, located at $w$. Define $\eta \wedge \eta'$ as the vertex-wise minimum of $\eta$ and $\eta'$, so that $\eta \wedge \eta' = \eta - \chi_u = \eta' - \chi_v$. We will use $\eta \to \eta'$ flows that pass only through configurations of the form $\eta \wedge \eta' + \chi_w$; that is, the flow-carrying paths send the "extra particle" at $u$ to $v$ along a path in $G_1$. To specify such a flow, it is enough to specify a flow $g_{uv}$ from $u$ to $v$ in $G_1$. We will use the one that simply spreads flow uniformly over all shortest paths from $u$ to $v$. Note that if $(\zeta_1, \zeta_2)$ is an edge in a shortest path from $\eta$ to $\eta'$, then $\zeta_1 \wedge \zeta_2 = \eta \wedge \eta'$.

Fix an edge $e = (\zeta_1, \zeta_2) \in \mathcal{E}_1$. We must bound the total flow through $e$. Let $\zeta = \zeta_1 \wedge \zeta_2$. Then all pairs $(\eta, \eta')$ such that some $\eta \to \eta'$ flow passes through $e$ must satisfy $\eta \wedge \eta' = \zeta$. For vertices $w$, define $\zeta^w = \zeta + \chi_w$. Then

$$f(e) = \sum_{u,v} f_{\zeta^u, \zeta^v}(e)$$
$$= \sum_{u,v} g_{uv}(e).$$

But Lemma 2 below says that $\sum_{u,v} g_{uv}(e) \leq L(L^d - 1)$. It follows that $\mathcal{C}(f) \leq L$. Finally, note that the flow $f$ also satisfies $\mathcal{L}(f) \leq dL$. Hence, Theorem 1 implies that $\tau_1 \leq 2d^2 L^2 \tau_2$.

The following lemma was used in the estimate for $\mathcal{C}(f)$:

LEMMA 2. *For vertices $u, v \in \mathbf{Z}^d / L\mathbf{Z}^d$, let $g_{uv}$ be the flow that spreads flow evenly among all shortest paths from $u$ to $v$. Then, for any edge $e$,*

$$\sum_{u,v} g_{uv}(e) \leq L(L^d - 1).$$

PROOF. Define $g = \sum_{u,v} g_{uv}$. By symmetry, the quantity $g(e)$ does not depend on $e$. But the maximum length of any shortest path is less than $dL$, the number of ordered pairs of distinct vertices is $L^d(L^d - 1)$ and the total number of edges is $dL^d$. Hence $e$ carries at most $\frac{dL \times (L^d)(L^d - 1)}{dL^d} = L(L^d - 1)$ units of flow. □

**3. The ZRP on the complete graph.** In the previous section we showed that $\tau_1 \leq 2d^2 L^2 \tau_2$. Thus, to prove an $O((\rho + 1)^2 L^2)$ bound for $\tau_1$, we need only show that the relaxation time $\tau_2$ for the ZRP on the complete graph can be bounded by $O((\rho + 1)^2)$, uniformly in the number of vertices.

THEOREM 3. *Let $K_n$ denote the complete graph on $n$ vertices. Fix $\rho > 0$ and let $\tau$ be the relaxation time for the ZRP on $K_n$ with a density $\rho$ of particles. Then $\tau \leq C(\rho + 1)^2$ for a universal constant $C$.*

PROOF. We will estimate $\tau$ using coupling; that is, we construct a process $\langle (\eta_t, \eta'_t), t \geq 0 \rangle$ such that $\langle \eta_t, t \geq 0 \rangle$ and $\langle \eta'_t, t \geq 0 \rangle$ are each ZRP's and the distribution of $\eta'_0$ is uniform. Let $T = \inf\{t : \eta_t = \eta'_t\}$ be the coupling time. It is a standard fact that

(2) $$\|K^t(\eta, \cdot) - \mathcal{U}\| \leq \mathbf{P}_\eta(T > t),$$

where we write $\mathbf{P}_\eta(\cdot) := \mathbf{P}(\cdot \mid \eta_0 = \eta)$. The following lemma relates the relaxation time $\tau$ to the tail of the distribution of $T$:

LEMMA 4. *Fix $\gamma > 0$ and suppose that $\mathbf{E}(e^{\gamma T}) < \infty$. Then $\tau \leq 1/\gamma$.*

PROOF. If $\mathbf{E}(e^{\gamma T}) < \infty$, then $\lim_{t \to \infty} e^{\gamma t} \mathbf{P}(T > t) = 0$. Hence $\mathbf{P}(T > t) < e^{-\gamma t}$ for sufficiently large $t$. Combining this with equation (2) gives $\tau \leq 1/\gamma$. □

*Zero range process with ranked particles.* We will actually study a process with labeled particles. Suppose that the particles are labeled $1, \ldots, r$ and each particle $x$ has a rank $R_t(x) \in \{1, \ldots, r\}$ at time $t$. Every particle has a unique rank, so $R_t(\cdot)$ is a permutation. We assume that when vertex $u$ rings, the highest ranking particle at $u$ (i.e., the particle whose rank is the smallest number) is expelled. Let $v_t(i)$ denote the location of the particle of rank $i$ at time $t$. Then $\eta_t(u) = \sum_{i=1}^{r} \mathbf{1}(v_t(i) = u)$, where we write $\mathbf{1}(A)$ for the indicator of the event $A$. For positive integers $j$, we define $\eta_t(u, j) = \sum_{i=1}^{j} \mathbf{1}(v_t(i) = u)$. We will use the notation $\eta_t(\cdot, \cdot)$ to denote a ZRP with ranked particles.

We allow the rankings to vary in time. However, we assume that the vertices ring independently of the rankings. Thus, if the particles of rank $1, \ldots, j$ maintain a constant rank [or even if the set $\{x : 1 \leq R(x) \leq j\}$ does not change], then the process $\langle \eta_t(\cdot, j) : t \geq 0 \rangle$ behaves like the ZRP with $j$ particles.

*The coupling.* We now describe the coupling. The coupling will consist of $r$ stages, going from Stage 0 to Stage $r - 1$. At the end of Stage $j$, configurations $\eta_t$ and $\eta'_t$ will satisfy $\eta_t(\cdot, j+1) = \eta'_t(\cdot, j+1)$. Thus, we will have $\eta_t = \eta'_t$ at the end of Stage $r - 1$. We will now describe Stage $r - 1$. (Earlier stages will be similar.) The stage consists of two steps:

- *Step* 1. First, we rank the particles so that the location of each particle whose rank is less than $r$ is the same in $\eta$ as it is in $\eta'$. Throughout Step 1, we couple the processes so that corresponding vertices ring at the same time, and when particles are expelled, they choose the same destination. This ensures that the locations of the particles of rank $1, 2, \ldots, r-1$ remain matched. We run the first step until
$$\eta_t(v_t(r), r-1) = \eta'_t(v'_t(r), r-1),$$
that is, the number of particles in the same location as the particle of rank $r$ is the same in $\eta$ as it is in $\eta'$. Let $a = v_t(r)$ and $b = v'_t(r)$. Note that

(3) $$\eta_t(a) = \eta'_t(b)$$

and

(4) $$\eta_t(b) = \eta'_t(a)$$

at the completion of Step 1.

- *Step* 2. In Step 2, we couple so that vertex $a$ in $\eta$ plays the role of vertex $b$ in $\eta'$, and vice versa. This is accomplished as follows. First, the particles located at vertices $a$ and $b$ are given new rankings so that the particles of vertex $a$ (resp. $b$) in $\eta$ correspond to the particles of vertex $b$ (resp. $a$) in $\eta'$. Second, we couple so that vertex $a$ (resp. $b$) rings in the process $\eta_t$ if and only if $b$ (resp. $a$) rings in the process $\eta'_t$. If an expelled particle moves to $a$ (resp. $b$) in $\eta_t$, then the corresponding particle moves to $b$ (resp. $a$) in $\eta'_t$, and so on. This ensures that (3) and (4) persist throughout the duration of Step 2. We run Step 2 until $\eta_t(a) = \eta_t(b)$. Note that this implies $\eta_t(\cdot, r) = \eta'_t(\cdot, r)$, and we are done.

Stage $j$ for $1 \leq j < r-1$ can be described similarly, except that in Stage $j$, only the trajectories of particles $1, \ldots, j+1$ are coupled. Particles $j+2, \ldots, r$ are allowed to move in an arbitrary way. Note that since particles $j+2, \ldots, r$ move only when there are no higher ranking particles in the same location, they do not interfere with the trajectories of particles $1, \ldots, j+1$, so Stage $j$ can proceed in similar fashion to Stage $r-1$. Stage $j$ finishes when $\eta(\cdot, j+1) = \eta'(\cdot, j+1)$.

Let $\zeta_t$ be an encoding of the coupled pair of labeled particle processes and the rankings at time $t$, and let $\text{stage}_t$ and $\text{step}_t$ denote the stage and step, respectively, in progress at time $t$. Then the process $\langle (\zeta_t, \text{stage}_t, \text{step}_t) : t \geq 0 \rangle$ is a Markov chain. We wish to show that for a universal constant $C > 0$, the time $\tau$ to complete stages 1 through $r-1$ satisfies $\mathbf{E}(e^{\gamma \tau}) < \infty$ for $\gamma \geq 1/C(\rho+1)^2$. The following lemma says that it is enough to show that for any particular step of any stage, the time $T$ spent there satisfies $\mathbf{E}(e^{\gamma T}) < \infty$ for $\gamma \geq 1/C(\rho+1)^2$:

LEMMA 5. *Let $X_t$ be a finite-state Markov chain with $k$ transient classes $S_1, S_2, \ldots, S_k$, and suppose that there are no transitions from $S_j$ to $S_{j-1}$ for $j \geq 2$. Define $T_j = \inf\{t > 0 : X_t \notin S_1 \cup \cdots \cup S_j\}$. Fix $\gamma > 0$ and suppose that for all $j$ and $x \in S_j$, we have $\mathbf{E}_x(e^{\gamma T_j}) < \infty$. Then $\mathbf{E}(e^{\gamma T_k}) < \infty$, irrespective of the distribution of $X_0$.*

PROOF. Define $T_0 = 0$, and for $j \in \{1, \ldots, k\}$, let $W_j = T_j - T_{j-1}$. Then $T_k = \sum_{j=1}^k W_j$. So $\mathbf{E}(e^{\gamma T_k}) = \mathbf{E}(e^{\gamma(W_1 + \cdots + W_k)})$. Let $\alpha_j = \max_{x \in S_j} \mathbf{E}_x(e^{\gamma T_j})$. We show by mathematical induction that $\mathbf{E}(e^{\gamma(W_1 + \cdots + W_l)}) \leq \prod_{j=1}^l \alpha_j$ for all $l \leq k$. The base case $l = 1$ follows from the definition of $\alpha_1$. Now fix $l < k$ and suppose that $\mathbf{E}(e^{\gamma(W_1 + \cdots + W_l)}) \leq \prod_{j=1}^l \alpha_j$. Let $Y$ be the first state that $X_t$ visits after leaving $S_l$. Then,

$$\mathbf{E}(e^{\gamma(W_1 + \cdots + W_{l+1})}) = \mathbf{E}(\mathbf{E}(e^{\gamma(W_1 + \cdots + W_{l+1})} | W_1, \ldots, W_l, Y))$$
$$= \mathbf{E}(e^{\gamma(W_1 + \cdots + W_l)} \mathbf{E}(e^{\gamma W_{l+1}} | W_1, \ldots, W_l, Y))$$

$$\leq \mathbf{E}\left(e^{\gamma(W_1+\cdots+W_l)} \max_{x \in S_{l+1}} \mathbf{E}_x(e^{\gamma T_{l+1}})\right)$$

$$= \alpha_{j+1}\mathbf{E}(e^{\gamma(W_1+\cdots+W_l)}),$$

and the result follows. $\square$

To complete the proof of Theorem 3, we need to show that the time $T$ spent performing a particular step satisfies $\mathbf{E}(e^{\gamma T}) < \infty$ for $\gamma \geq 1/C(\rho+1)^2$. Lemmas 9 and 10 below do this for steps of type 2 and 1, respectively. $\square$

We will need the following lemma:

LEMMA 6. *Fix $v \in V$ and suppose that $\eta_0(v) \leq 2(\rho+1)$. Let $Z = |\{t < (\rho+1)^2 : \eta_v(t) = 0\}|$. Then, for universal constants $C > 0$ and $M > 0$, we have $\mathbf{E}(Z \wedge M(\rho+1)) \geq C(\rho+1)$.*

PROOF. Note that the ZRP on $K_n$ can be described in the following way. For each ordered pair of distinct vertices $(u, v)$ at rate $1/(n-1)$, the process attempts to move a particle from $u$ to $v$; that is, the move is made unless $u$ is empty. Thus, for any vertex $v$, the total rate at which attempts to increase $\eta_t(v)$ are made is 1, and this is also the rate of attempts to decrease $\eta_t(v)$. Furthermore, when $\eta_t(v) > 0$ it is always possible to decrease $\eta_t(v)$, so $\eta_t(v) \to \eta_t(v) - 1$ at rate 1 in this case. Let $i = \eta_0(v)$ and let $X_t$ be the continuous-time simple symmetric random walk on the integers, started at $i$, that moves in each direction at rate 1.

One can easily couple the processes $\eta_t(v)$ and $X_t$ so that $\eta_t(v) = 0$ whenever $X_t = 0$. Thus, it is enough to show that for constants $C > 0$ and $M > 0$, we have $\mathbf{E}(Y \wedge M(\rho+1)) \geq C(\rho+1)$, where $Y = |\{t < (\rho+1)^2 : X_t = 0\}|$.

We will first show that $\mathbf{E}(Y) \geq D(\rho+1)$ for a universal constant $D$ by showing that if $t$ satisfies $\frac{1}{2}(\rho+1)^2 \leq t \leq (\rho+1)^2$, then $\mathbf{P}(X_t = 0) \geq 2D(\rho+1)^{-1}$.

Note that $X_t - i$ is distributed as the difference of two independent Poisson($t$) random variables. Hence, $X_t - i$ is approximately normal with mean 0 and variance $2t$ when $t$ is large. More precisely, the local limit theorem implies that there is a universal constant $D$ such that for all $t \geq \frac{1}{2}$ and $j$ with $|j| \leq \sqrt{8t}$, we have

$$\mathbf{P}(X_t - i = j) \geq \frac{2D}{\sqrt{t}}. \tag{5}$$

Recall that $i \leq 2(\rho+1)$ and note that if $t \in [\frac{1}{2}(\rho+1)^2, (\rho+1)^2]$, then $\sqrt{t} \leq (\rho+1) \leq \sqrt{2t}$. Hence $i \leq \sqrt{8t}$. Thus, we can substitute $j = -i$ into (5) and obtain

$$\mathbf{P}(X_t = 0) \geq 2D/\sqrt{t} \geq 2D/(\rho+1).$$

Next, we will show that $Y$ is unlikely to be a very large multiple of $(\rho+1)$. Let $T_1 = \inf\{t : X_t = 0\}$ and for $k \geq 2$, let $T_k = \inf\{t > T_{k-1}+1 : X_t = 0\}$. For $k \geq 1$, let $A_k$ be the event that $X_t \neq 0$ for all $t \in [T_k + 1, T_k + (\rho+1)^2]$.

Proposition A.2 in the Appendix gives

$$\mathbf{P}_0(X_t \neq 0 \text{ for all } t \in [1, (\rho+1)^2]) \geq b(\rho+1)^{-1},$$

for a universal constant $b$. Combining this with the strong Markov property gives

$$\mathbf{P}(A_k \,|\, A_1^c \cdots A_{k-1}^c) \geq b(\rho+1)^{-1},$$

and hence $\mathbf{P}(A_1^c \cdots A_k^c) \leq (1 - b(\rho+1)^{-1})^k$.

Recall that $Y$ is the amount of time $X_t$ spends at 0 up to time $(\rho+1)^2$. Note that $Y \leq N$, where $N = \inf\{k : T_k \geq (\rho+1)^2\}$. Hence, for every nonnegative integer $k$, we have

$$\mathbf{P}(Y > k) \leq \mathbf{P}(N > k)$$
$$\leq \mathbf{P}(A_1^c \cdots A_k^c)$$
$$\leq (1 - b(\rho+1)^{-1})^k.$$

Let $M = b^{-1} \log(2/Db)$. Note that

$$Y = (Y \wedge M(\rho+1)) + W,$$

where $W = (Y - M(\rho+1))\mathbf{1}(Y > M(\rho+1))$. Since $\mathbf{E}(Y) \geq D(\rho+1)$, it follows that

(6) $\qquad \mathbf{E}(Y \wedge M(\rho+1)) \geq D(\rho+1) - \mathbf{E}(W).$

But

$$\mathbf{E}(W) \leq \sum_{k \geq 0} \mathbf{P}(W > k)$$
$$= \sum_{k \geq 0} \mathbf{P}(Y > M(\rho+1) + k)$$
$$\leq \sum_{k \geq 0} (1 - b(\rho+1)^{-1})^{M(\rho+1)+k}.$$

By summing the geometric series and then using the fact that $1 + u \leq e^u$ for all $u$, it is easily shown that $\mathbf{E}(W) \leq e^{-bM} b^{-1}(\rho+1) = \frac{D}{2}(\rho+1)$. Combining this with (6) gives $\mathbf{E}(Y \wedge M(\rho+1)) \geq \frac{D}{2}(\rho+1)$. □

The following lemma is a consequence of Lemma 6. It says that when $t$ is large, the average fraction of time the vertices spend empty is $\Omega(1/(\rho+1))$, with high probability. For $v \in V$, define

$$Z_t(v) = |\{s < t : \eta_s(v) = 0\}|$$

and let $\overline{Z}_t = \frac{1}{n} \sum_{v \in V} Z_t(v)$.

LEMMA 7. *There exists a universal constant $\alpha > 0$ such that for every $t > (\rho+1)^2$, we have*

$$\mathbf{P}\left(\overline{Z}_t < \frac{Ct}{4(\rho+1)}\right) \leq e^{-\alpha t/(\rho+1)^2},$$

*where $C$ is the constant appearing in Lemma 6.*

PROOF. Let $M$ be the constant appearing in Lemma 6. For $k \in \{1, 2, \ldots\}$, let

$$X_k(v) = |\{(k-1)(\rho+1)^2 \leq s \leq k(\rho+1)^2 : \eta_s(v) = 0\}| \wedge M(\rho+1)$$

and let $\overline{X}_k = \frac{1}{n}\sum_{v \in V} X_k(v)$. Markov's inequality implies that at least half of the vertices $v \in V$ satisfy $\eta_0(v) \leq 2(\rho+1)$. For such vertices $v$, Lemma 6 says that $\mathbf{E}(X_1(v)) \geq C(\rho+1)$. Hence $\mathbf{E}(\overline{X}_1) \geq \frac{1}{2}C(\rho+1)$. A similar argument applies to $\overline{X}_k$, for every positive integer $k$. It follows that if we define $S_k = \sum_{i=1}^{k}(\overline{X}_i - \frac{1}{2}C(\rho+1))$, then $\langle S_k : k \geq 0 \rangle$ is a submartingale. Recall Azuma's inequality (see Corollary 6.9 on page 166 and Section (c) on page 168 of [8]): if $M_k$ is a submartingale with $M_0 = 0$ and $|M_{k+1} - M_k| \leq B$ for all $k$, then for every $A > 0$ and $k \geq 1$, we have

$$\mathbf{P}\left(\min_{j \leq k} M_j \leq -A\right) \leq e^{-bA^2/kB^2},$$

where $b > 0$ is a universal constant. (In [8], Azuma's inequality is stated for martingales, but since any submartingale can be written as a martingale plus a nonnegative, predictable sequence, clearly the result must hold for submartingales as well.) Note that $\{S_k\}$ satisfies $|S_{k+1} - S_k| \leq (M + \frac{1}{2}C)(\rho+1)$. Hence, Azuma's inequality implies that

$$\mathbf{P}(S_k \leq -A) \leq e^{-b'A^2/k(\rho+1)^2},$$

for a universal constant $b' > 0$.

Hence,

$$\mathbf{P}\left(\sum_{i=1}^{k} \overline{X}_i \leq \tfrac{1}{4}kC(\rho+1)\right) = \mathbf{P}(S_k \leq -\tfrac{1}{4}kC(\rho+1)) \leq e^{-b'C^2k/16}.$$

Note that $\overline{Z}_{k(\rho+1)^2} \geq \sum_{i=1}^{k} \overline{X}_i$. Thus, for $t = k(\rho+1)^2$ and $\alpha = b'C^2/16$, we have

(7) $\quad \mathbf{P}\left(\overline{Z}_t \leq \dfrac{Ct}{4(\rho+1)}\right) \leq \mathbf{P}\left(\sum_{i=1}^{k} \overline{X}_i \leq \dfrac{Ct}{4(\rho+1)}\right) \leq e^{-\alpha t/(\rho+1)^2}.$

Thus the lemma holds for $t$ of the form $k(\rho+1)^2$. But since $\overline{Z}_t$ is nondecreasing in $t$, we can make (7) true for all $t \geq (\rho+1)^2$ by incorporating an extra factor of $\frac{1}{2}$ into $\alpha$. □

We also use the following proposition:

PROPOSITION 8. *There is a universal constant $A > 0$ such that if $X$ and $Y$ are independent* Poisson$(\lambda)$ *random variables for sufficiently large $\lambda$, then:*

(a) $\mathbf{P}(|X - \lambda| \geq \lambda/2) \leq e^{-A\lambda}$;
(b) $\mathbf{P}(X - Y \geq \alpha\lambda) \leq e^{-A\alpha^2\lambda}$ *for every $\alpha \in [0, 1]$.*

PROOF. We prove (b) first. The moment generating function of $X$ is
$$\phi_X(\theta) \equiv \mathbf{E}(e^{\theta X}) = \exp(\lambda(e^\theta - 1)).$$
It follows that the moment generating function of $X - Y$ is
$$\phi_{X-Y}(\theta) = \phi_X(\theta)\phi_X(-\theta)$$
$$= \exp(\lambda(e^\theta - 1) + \lambda(e^{-\theta} - 1))$$
$$= \exp(\lambda(\theta^2 + o(\theta^2))),$$
since $e^u = 1 + u + \frac{1}{2}u^2 + o(u^2)$. Hence,
$$\mathbf{P}(X - Y \geq \alpha\lambda) = \mathbf{P}(e^{\theta(X-Y)} \geq e^{\theta\alpha\lambda})$$
$$\leq e^{-\theta\alpha\lambda}\phi_{X-Y}(\theta)$$
$$= e^{-\theta\alpha\lambda}e^{\lambda(\theta^2 + o(\theta^2))},$$
by Markov's inequality. Taking logarithms, we get
$$\log(\mathbf{P}(X - Y \geq \alpha\lambda)) \leq \lambda(\theta^2 - \theta\alpha + o(\theta^2)).$$
When $\theta = \alpha/2$, this becomes
$$\lambda\left(\frac{\alpha^2}{4} - \frac{\alpha^2}{2} + o(\alpha^2)\right) = \lambda\left(-\frac{1}{4}\alpha^2 + o(\alpha^2)\right).$$
Hence $\mathbf{P}(X - Y \geq \alpha\lambda) \leq e^{-A\alpha^2\lambda}$ for a universal constant $A$ and (b) is proved.

For (a), note that
$$\phi_{X-\lambda}(\theta) = \exp(\lambda(e^\theta - \theta - 1)); \qquad \phi_{\lambda-X}(\theta) = \exp(\lambda(e^{-\theta} + \theta - 1)).$$
Using the inequality $e^{|u|} \leq e^u + e^{-u}$ and setting $\theta = \frac{1}{2}$ gives
$$\phi_{|X-\lambda|}(\tfrac{1}{2}) \leq \exp(\lambda(e^{1/2} - \tfrac{1}{2} - 1)) + \exp(\lambda(e^{-1/2} + \tfrac{1}{2} - 1))$$
$$\leq 2e^{\lambda/10},$$
since $e^u - u - 1 \leq 1/10$ whenever $u \in \{-\frac{1}{2}, \frac{1}{2}\}$. Applying Markov's inequality to the random variable $e^{|X-\lambda|/2}$ gives
$$\mathbf{P}(|X - \lambda| \geq \lambda/2) \leq e^{-\lambda/4} \cdot 2e^{\lambda/10} = 2e^{-3\lambda/20},$$
which is at most $e^{-\lambda/10}$ when $\lambda$ is sufficiently large. □

We will now use Lemma 7 to prove Lemmas 9 and 10. It turns out that steps of type 2 are easier to analyze, so we start with the following:

LEMMA 9. *Fix $j \in \{0, \ldots, r-1\}$ and let $T$ be the time spent performing Step 2 of Stage $j$. Then $\mathbf{E}(e^{\gamma T}) < \infty$ for $\gamma \geq D(\rho+1)^{-2}$, where $D$ is a universal constant.*

PROOF. We may assume that $r = j+1$, since if there are particles that rank lower than $j+1$, they do not affect the behavior of the particles ranking $1, \ldots, j+1$ and they increase the value of $\rho$.

Let vertices $a$ and $b$ be as defined in the description of the coupling and suppose that $\eta_t(a) \geq \eta_t(b)$ at the beginning of Step 2. Step 2 finishes when $\eta_t(a) = \eta_t(b)$, so Step 2 has been completed by the time that $\eta_t(a) = 0$. Thus, it is enough to show that for any vertex $v$, if we define $\tau = \inf\{t : \eta_t(v) = 0\}$, then $\mathbf{E}(e^{\gamma \tau}) < \infty$ for $\gamma \geq D(\rho+1)^{-2}$.

Fix a vertex $v$, and for $t \geq 0$, let $X_t$ be the number of attempts to increase the occupancy of $v$ minus the number of attempts to decrease the occupancy of $v$, up to time $t$. Let $B_t$ be the number of attempts to increase the occupancy of $v$, up to time $t$, that fail because the potential source vertex is empty. Note that before time $\tau$, all attempts to decrease the occupancy of $v$ will succeed. Therefore,

$$\eta_t(v) = \eta_0(v) + X_t - B_t,$$

up to time $\tau$. It follows that $\mathbf{P}(\tau > t) \leq \mathbf{P}(X_t - B_t > -\eta_0(v))$.

We assume that even when a vertex is empty, it attempts to move a particle to $v$ at rate $\frac{1}{n-1} > \frac{1}{n}$. [We will use the convention that when $v$ is empty, it attempts to move a particle to itself, and thereby increases $B_t$, at rate $1/(n-1)$.] It follows that, given $\overline{Z}_t$, the conditional distribution of $B_t$ stochastically dominates a Poisson($\overline{Z}_t$) random variable. (Recall that the sum of the amounts of time the vertices have been empty up to time $t$ is $Z_t = n\overline{Z}_t$.) Thus, for $t$ sufficiently large, we have

(8)
$$\mathbf{P}\left(B_t \leq \frac{Ct}{8(\rho+1)}\right) \leq \mathbf{P}\left(B_t \leq \frac{Ct}{8(\rho+1)} \Big| \overline{Z}_t \geq \frac{Ct}{4(\rho+1)}\right)$$
$$+ \mathbf{P}\left(\overline{Z}_t < \frac{Ct}{4(\rho+1)}\right)$$

(9)
$$\leq e^{-AtC/4(\rho+1)} + e^{-\alpha t/(\rho+1)^2},$$

where the second inequality uses part (a) of Proposition 8 (for the first term) and Lemma 7 (for the second term). It follows that for sufficiently large $t$, we have

(10)
$$\mathbf{P}\left(B_t \leq \frac{Ct}{8(\rho+1)}\right) \leq \exp\{-ct/(\rho+1)^2\},$$

for a universal constant $c$. Since $X_t$ is the difference of two independent Poisson($t$) random variables, we can apply part (b) of Proposition 8 to $X_t$

and get

(11) $$\mathbf{P}\left(X_t \geq \frac{Ct}{16(\rho+1)}\right) \leq \exp\{-AC^2 t/256(\rho+1)^2\}.$$

Equations (10) and (11) together imply that for some constant $\delta > 0$, we have

$$\mathbf{P}\left(X_t - B_t > -\frac{Ct}{16(\rho+1)}\right) \leq \exp\{-\delta t/(\rho+1)^2\},$$

for sufficiently large $t$, and so

$$\mathbf{P}(\tau > t) \leq \mathbf{P}(X_t - B_t > -\eta_0(v))$$
$$\leq \exp\{-\delta t/(\rho+1)^2\},$$

for sufficiently large $t$ [e.g., for $t$ such that $Ct/16(\rho+1) > r$]. It follows that $\mathbf{E}(e^{\delta \tau/2(\rho+1)^2}) < \infty$. $\square$

Next, we consider steps of type 1.

LEMMA 10. *Fix $j \in \{0,\ldots,r-1\}$ and let $T$ be the time spent performing Step 1 of Stage $j$ of the coupling. Then $\mathbf{E}(e^{\gamma T}) < \infty$ for $\gamma \geq D(\rho+1)^{-2}$, where $D$ is a universal constant.*

PROOF. Again we may assume that $r = j+1$. In order to analyze Step 1, we will embed a coupling into a *single* process as follows. Consider a zero range process $\eta_t(\cdot,\cdot)$ with $j+2$ ranked particles, and for $t \geq 0$, define $\xi_t$ and $\xi'_t$ by

$$\xi_t = \eta_t(\cdot,j) + \chi_{v_t(j+1)}, \qquad \xi'_t = \eta_t(\cdot,j) + \chi_{v_t(j+2)}.$$

So $\xi_t$ and $\xi'_t$ are obtained from $\eta_t(\cdot,\cdot)$ by deleting the particles of rank $j+2$ and $j+1$, respectively. Let $\eta_0(\cdot,\cdot)$ be defined so that $(\xi_0, \xi'_0)$ reflects the state of the coupling at the beginning of Step 1 of Stage $j$. Then the time $T$ to complete Step 1 has the same distribution as $W$, where $W = \inf\{t : \eta_t(v_t(j+1), j+2) = \eta_t(v_t(j+2), j+2)\}$, because up to time $W$, the process $(\xi_t, \xi'_t)$ behaves exactly like the coupling during Step 1. [Note that the particle of rank $j+1$ cannot interfere with the particle of rank $j+2$ because $W \leq \inf\{t : v_t(j+1) = v_t(j+2)\}$.]

We wish to show that $\mathbf{E}(e^{\gamma W})$ is finite for $\gamma = D(\rho+1)^{-2}$. Consider the process $\langle \zeta_t : t \geq 0 \rangle$, where for $t \geq 0$, we define $\zeta_t = (\eta_t(\cdot), v_t(j+1), v_t(j+2))$. [Here, we write $\eta_t(\cdot)$ for the function that outputs the number of particles at each site, but contains no information about their rankings. So $\zeta_t$ encodes only the number of particles at each vertex at time $t$ and the positions of the two lowest ranking particles.] The particles of rank $1, 2, \ldots, j$ will be called *high ranking*; note that in the process $\zeta_t$, the high ranking particles are indistinguishable. Furthermore, $\{\zeta_t : t \geq 0\}$ is a Markov chain.

It will be convenient to study a modified version of the process obtained by combining all of the states of the form $(\eta, u, u)$ into a single state, which we denote by $b$. (Note that since the chain does not reach $b$ before time $W$, making such a modification does not affect the distribution of $W$.) We shall also choose transition rates $q_{b,x}$ out of $b$ that ensure a "nice" stationary distribution. (Again, this does not affect $W$.) For all $x = (\eta, s, t)$ with $s \neq t$, let

$$
(12) \qquad q_{b,x} = \begin{cases} \dfrac{1}{n-1}, & \text{if } \eta(s) = \eta(t) = 1, \\ \dfrac{2}{n-1}, & \text{otherwise.} \end{cases}
$$

Recall that if $\pi$ is a positive function on the state space that solves the balance equations

$$
(13) \qquad \sum_y \pi(y) q_{y,x} = \sum_y \pi(x) q_{x,y}
$$

for every $x$, then $\pi$ can be normalized to be the stationary distribution. [In fact, it is enough to verify (13) for all but one of the states $x$.] We claim that if $\pi$ is defined by $\pi(b) = 1$ and

$$
(14) \qquad \pi(\eta, s, t) = \eta(s)\eta(t)
$$

when $s \neq t$, then $\pi$ solves (13). For all $x$, define $Q_{\text{in}}(x) = \sum_y \pi(y) q_{y,x}$ and $Q_{\text{out}}(x) = \sum_y \pi(x) q_{x,y}$, so that (13) is equivalent to $Q_{\text{in}}(x) = Q_{\text{out}}(x)$ for all $x$. We will verify this for all $x \neq b$. For configurations $x$ and vertices $u$ and $v$, define $Q_{\text{out}}(x; \{u, v\}) = \sum_y \pi(x) q_{x,y}$, where the sum is over configurations $y$ obtained from $x$ by either moving a particle from $u$ to $v$ or vice versa, with a similar definition for $Q_{\text{in}}(x; \{u, v\})$. Note that

$$
(15) \qquad Q_{\text{out}}(x) = \sum_{\{u,v\}} Q_{\text{out}}(x; \{u, v\}),
$$

$$
(16) \qquad Q_{\text{in}}(x) = \sum_{\{u,v\}} Q_{\text{in}}(x; \{u, v\}) + q_{b,x}.
$$

Fix $x = (\eta, s, t)$ with $s \neq t$. In order to calculate the values of $Q_{\text{out}}(x, \{u, v\})$, we shall consider three cases separately:

*Case* $\{u, v\} \cap \{s, t\} = \varnothing$. In this case, it is easy to verify that

$$
Q_{\text{out}}(x, \{u, v\}) = Q_{\text{in}}(x, \{u, v\})
$$
$$
= (\delta(u) + \delta(v))\frac{\eta(s)\eta(t)}{n-1},
$$

where for vertices $w$, we define
$$\delta(w) = \begin{cases} 1, & \text{if } \eta(w) \geq 1, \\ 0, & \text{otherwise.} \end{cases}$$

*Case $u = s, v \neq t$.* In this case, we have
$$Q_{\text{out}}(x, \{u, v\}) = (1 + \delta(v)) \frac{\eta(s)\eta(t)}{n-1}.$$

If $\eta(v) = 0$, then
$$Q_{\text{in}}(x, \{u, v\}) = \frac{(\eta(s) - 1)\eta(t)}{n-1} + \frac{\eta(t)}{n-1}$$
$$= \frac{\eta(s)\eta(t)}{n-1}.$$

If $\eta(v) \geq 1$, then
$$Q_{\text{in}}(x, \{u, v\}) = \frac{(\eta(s) - 1)\eta(t)}{n-1} + \frac{(\eta(s) + 1)\eta(t)}{n-1}$$
$$= \frac{2\eta(s)\eta(t)}{n-1}.$$

In both cases, $Q_{\text{out}}(x, \{u, v\}) = Q_{\text{in}}(x, \{u, v\})$.

*Case $u = s, v = t$.* This is the only case where $Q_{\text{out}}(x, \{u, v\}) \neq Q_{\text{in}}(x, \{u, v\})$. We have
$$Q_{\text{out}}(x, \{u, v\}) = \begin{cases} \dfrac{1}{n-1}, & \text{if } \eta(s) = \eta(t) = 1, \\ \dfrac{2\eta(s)\eta(t)}{n-1}, & \text{otherwise.} \end{cases}$$

If $\eta(s) \geq 2$ and $\eta(t) \geq 2$, then
$$Q_{\text{in}}(x, \{u, v\}) = \frac{(\eta(s) - 1)(\eta(t) + 1)}{n-1} + \frac{(\eta(s) + 1)(\eta(t) - 1)}{n-1}$$
$$= \frac{2\eta(s)\eta(t) - 2}{n-1}.$$

If $\eta(s) = 1$ and $\eta(t) \geq 2$, then
$$Q_{\text{in}}(x, \{u, v\}) = \frac{2(\eta(t) - 1)}{n-1}$$

and if $\eta(s) = \eta(t) = 1$, then $Q_{\text{in}}(x, \{u, v\}) = 0$. It follows that
$$Q_{\text{out}}(x, \{u, v\}) - Q_{\text{in}}(x, \{u, v\}) = \begin{cases} \dfrac{1}{n-1}, & \text{if } \eta(s) = \eta(t) = 1, \\ \dfrac{2}{n-1}, & \text{otherwise.} \end{cases}$$

Putting all of this together with equations (12), (15) and (16) verifies the balance equations.

Since the stationary distribution puts positive mass on every state, it is enough to show that $\mathbf{E}(e^{\gamma W})$ is finite if the chain starts in its stationary distribution.

Let $q_{i,j}$ be the rate at which $\zeta_t$ goes from $i$ to $j$. The time reversal $\tilde{\zeta}_t$ of $\zeta_t$ is the process that starts in distribution $\pi$ and has transition rates $\tilde{q}$ given by $\pi(i)q_{i,j} = \pi(j)\tilde{q}_{j,i}$. For events $A$, we will write $\tilde{\mathbf{P}}(A)$ to denote the probability of $A$ when we run $\tilde{\zeta}_t$ instead of $\zeta_t$, with a similar notation for expectation. Let $B$ be the set of states $(\eta, u, w)$ of $\zeta_t$ satisfying $\eta(u) = \eta(w)$, including the amalgamated state $b$. Then

$$\mathbf{P}(W > t') = \mathbf{P}(\zeta_t \notin B \text{ for all } t \in [0, t'])$$
$$= \tilde{\mathbf{P}}(\tilde{\zeta}_t \notin B \text{ for all } t \in [0, t']) = \tilde{\mathbf{P}}(W > t').$$

Thus, it is enough to show that $\tilde{\mathbf{E}}(e^{\gamma W}) < \infty$. Furthermore, we may modify the process $\tilde{\zeta}_t$ so that all transitions to $b$ are suppressed, because this can only increase $W$ (since $b \in B$). If we make this modification, it is then straightforward to verify, using (14), that the transition rule for the resulting process can be described as follows. Suppose that the current state is $(\eta, u_1, u_2)$ with $u_1 \neq u_2$. Then, for every ordered pair of distinct vertices $(v, w)$ at rate $\left(\frac{\eta(w,j+2)+1}{(\eta(w,j)+1)}\right) \times \frac{1}{n-1}$, an attempt is made to move a particle from $v$ to $w$, where an attempt proceeds as follows. If $v$ is empty, or the transition would put the low ranking particles in the same location, then the attempt fails. Otherwise, a particle is chosen "uniformly" from $v$ [i.e., a high ranking particle is chosen with probability $\frac{\eta(v,j)}{\eta(v,j+2)}$ and a low ranking particle is chosen with the remaining probability] and moved to $w$ if the reversed move is a transition of $\zeta_t$ (i.e., the high ranking particles are always allowed to move and the low ranking particles are only allowed to move if the destination vertex is empty).

Note that $\tilde{\zeta}_t$ is somewhat similar to $\zeta_t$, except for the motions of the low ranking particles. In $\zeta_t$, the low ranking particles occasionally jump from an empty vertex to a nonempty one, whereas in $\tilde{\zeta}_t$ they occasionally jump from a nonempty vertex to an empty one. The high ranking particles move nearly the same way in both processes, except that in $\tilde{\zeta}_t$ they have a slight preference for moving to vertices that contain a low ranking particle.

The reader might wonder why we are studying the time reversal $\tilde{\zeta}_t$ instead of $\zeta_t$ itself. The reason is that in $\zeta_t$, since the occupancies of the vertices containing $j+1$ and $j+2$, respectively, can make big positive jumps, it is harder to rule out the possibility that it takes a long time for these occupancies to be the same.

Up to time $W$, the rate at which any nonempty vertex $v$ expels a particle is at least $1 - \frac{1}{\eta(v,j+2)}$, and the rate at which the process attempts to move a particle to $v$ is at most $1 + \frac{1}{\eta(v,j+2)}$. Let

$$(17) \qquad M_t = \max(\eta_t(v_t(j+1)), \eta_t(v_t(j+2))),$$

that is, $M_t$ is the maximum number of particles in a vertex with a low ranking particle. Before time $W$, we have $\eta_t(v_t(j+1)) \neq \eta_t(v_t(j+2))$. Hence, up to time $W$, the rate at which $M_t$ decreases when $M_t = x$ (i.e., the rate at which the process moves to a new state that decreases $M_t$) is at least $1 - \frac{1}{x}$, and the rate at which the process attempts to increase $M_t$ [i.e., the rate at which the process attempts to move a particle to the vertex that achieves the maximum in (17)] is at most $1 + \frac{1}{x}$. Let $B_t$ be the number of attempts to increase $M$ up to time $t$ that fail because the potential source vertex is empty. Note that when a vertex is empty, it attempts to move a particle to the vertex that achieves the maximum in (17) at rate at least $1/(n-1)$. It follows that, given $\overline{Z}_t$, the conditional distribution of $B_t$ stochastically dominates a Poisson($\overline{Z}_t$) random variable. Furthermore, since $\overline{Z}_t$ is in the $\sigma$-field generated by the amounts of time the process spends in each of its states up to time $t$, Lemma A.1 in the Appendix implies that the conclusion of Lemma 7 holds with $\tilde{\mathbf{P}}$ replacing $\mathbf{P}$. (Formally, Lemmas A.1 and 7 are not enough, because we altered the process $\tilde{\zeta}$ by constructing the amalgamated state $b$. However, since this change only takes effect after time $W$, it does not affect our conclusion here.)

Thus, as in the proof of Lemma 9 [see the equations leading up to equation (10)], we can apply part (a) of Proposition 8 to show that for $t$ sufficiently large, we have

$$(18) \qquad \tilde{\mathbf{P}}\left(B_t \leq \frac{Ct}{8(\rho+1)}\right) \leq \exp(-ct/(\rho+1)^2),$$

for a universal constant $c$.

Let $\tilde{\rho} = 64C^{-1}(\rho+1)$ and let $\alpha = \frac{1}{\tilde{\rho}(\tilde{\rho}+1)}$. For integers $j \geq 2$, let $w(j) = -\max(\frac{1}{j(j-1)}, \alpha)$. Define $\mathcal{W}(1) = 0$ and for $j \geq 2$, let $\mathcal{W}(j) = \sum_{i=2}^{j} w(j)$. When $M_t$ jumps downward from $j$ to $i$, where $i < j$, we define the *size* of the jump by $\sum_{k=i+1}^{j-1} w(k)$. Note that the size of a jump is zero unless $i < j-1$, and when the size of a jump is nonzero it is a negative number. Strictly negative jumps can only occur when a particle ranked either $j+1$ or $j+2$ moves to an empty vertex, or when the particle of rank $j+1$ moves to the vertex containing the particle of rank $j+2$.

Let $J_t$ be the sum of the sizes of all jumps up to time $t$ and let $U_t = \mathcal{W}(M_t) + J_t$. For $t \geq 0$, let $\mathcal{F}_t$ denote the $\sigma$-field generated by $\{(\tilde{\zeta}_s, B_s) : 0 \leq s \leq t\}$.

LEMMA 11. *Let $Y_t = U_t - \alpha B_t + \frac{4t\alpha}{\tilde{\rho}}$. Then $(Y_{t\wedge W}, \mathcal{F}_t)$ is a submartingale.*

PROOF. Suppose that $t < W$ and that the current configuration of $\tilde{\zeta}_t$ makes $M_t = x$. Let $v$ be the vertex that achieves the maximum in equation (17) and let $p$ be the the fraction of vertices in $V - \{v\}$ that are empty. Note that

$$\begin{cases} B_t \to B_t + 1, & \text{at rate at most } p\left(\frac{x+1}{x}\right), \\ U_t \to U_t + w(x+1), & \text{at rate at most } (1-p)\left(\frac{x+1}{x}\right), \\ U_t \to U_t - w(x), & \text{at rate at least } \frac{x-1}{x}. \end{cases}$$

Hence, $\lim_{\epsilon \to 0} \frac{1}{\epsilon}\tilde{\mathbf{E}}((U_{t+\epsilon} - \alpha B_{t+\epsilon}) - (U_t - \alpha B_t)|\mathcal{F}_t, t < W)$ is at least

$$-\left[p\left(\frac{x+1}{x}\right)\alpha + (1-p)\left(\frac{x+1}{x}\right)\max\left(\frac{1}{x(x+1)}, \alpha\right)\right]$$
$$+ \frac{x-1}{x}\max\left(\frac{1}{x(x-1)}, \alpha\right)$$
$$\geq -\left(\frac{x+1}{x}\right)\max\left(\frac{1}{x(x+1)}, \alpha\right) + \frac{x-1}{x}\max\left(\frac{1}{x(x-1)}, \alpha\right).$$

This is zero unless $\alpha > \frac{1}{x(x+1)}$, in which case it is at least

(19) $$-\left(\frac{x+1}{x}\right)\alpha + \left(\frac{x-1}{x}\right)\alpha = \frac{-2\alpha}{x}.$$

Recall that $\alpha = \frac{1}{\tilde{\rho}(\tilde{\rho}+1)}$. Thus, if $\alpha > \frac{1}{x(x+1)}$, then $x \geq \tilde{\rho}$, so $\frac{-2\alpha}{x} \geq \frac{-2\alpha}{\tilde{\rho}}$.

Define $\mathcal{N}_t = Y_{t \wedge W}$. The above calculation shows that

$$\lim_{\epsilon \to 0} \epsilon^{-1} \mathbf{E}(\mathcal{N}_{t+\epsilon} - \mathcal{N}_t | \mathcal{F}_t, t < W) \geq 4\alpha/\tilde{\rho} - 2\alpha/\tilde{\rho} = 2\alpha/\tilde{\rho}.$$

It is also clear that for some universal constant $\Delta > 0$, we have

$$\epsilon^{-1} \mathbf{E}(\mathcal{N}_{t+\epsilon} - \mathcal{N}_t | \mathcal{F}_t, t < W) \geq 0,$$

for every $t$, whenever $\epsilon < \Delta$. Thus, for every $t$, the discrete-time process $\langle \mathcal{N}_{t+k\delta} : k = 0, 1, \ldots \rangle$ is a submartingale whenever $\delta < \Delta$. It follows that for every $s > t$, we have $\mathbf{E}(\mathcal{N}_s | \mathcal{F}_t) \geq \mathcal{N}_t$, since $s$ can be written as $t + k\delta$ for some positive integer $k$, and $\langle \mathcal{N}_{t+k\delta} : k \geq 0 \rangle$ is a submartingale. It follows that $\mathcal{N}_t$ is a submartingale. $\square$

We will now briefly sketch how Lemma 11 will be used to prove Lemma 10. First, we will show that $\alpha B_t$ grows much faster than $\frac{4t\alpha}{\tilde{\rho}}$ as $t \to \infty$. But Lemma 11 and Azuma's inequality imply that $Y_{t \wedge W}$ is not likely to be much smaller than 0 when $t$ is large. Since $U_t \leq 0$, this means that $W$ is not likely to be very large.

More precisely, let $T_0 = 0$ and for $j \geq 1$, let $T_j = W \wedge \inf\{t > T_{j-1} : |Y_t - Y_{T_{j-1}}| \geq \frac{1}{2}\}$. For $k \geq 1$, let $M_k = Y_{T_k}$. Then $\{M_k : k \geq 0\}$ is a submartingale and $|M_{k+1} - M_k| \leq 1$ because the magnitudes of the jumps of $Y_t$ are at most $\frac{1}{2}$. Thus, Azuma's inequality implies that for any constant $A > 0$, we have

$$\tilde{\mathbf{P}}\left(\min_{j \leq k} M_j - M_0 < -A\right) \leq e^{-bA^2/k}$$

for a universal constant $b > 0$. Let $N_t = \inf\{k : T_k \geq t \wedge W\}$. The "drift" (viz., $\frac{4\alpha}{\tilde{\rho}}$) and "rate of jumping" [i.e., $\lim_{\epsilon \to 0} \epsilon^{-1} \tilde{\mathbf{P}}(Y_{t+\epsilon} \neq Y_t | \mathcal{F}_t)$] of $Y_t$ are uniformly bounded. Thus, for universal constants $h, \gamma > 0$, $N_t$ is stochastically dominated by $\gamma t + X$, for a Poisson($ht$) random variable $X$. Thus, Proposition 8 implies that $\tilde{\mathbf{P}}(N_t > \gamma t + 2ht) \leq e^{-Aht}$ for sufficiently large $t$. So if we define $\beta = \gamma + 2h$, then $\tilde{\mathbf{P}}(N_t > \beta t) \leq e^{-Aht}$ for sufficiently large $t$. Note that for any $\epsilon > 0$, we have for sufficiently large $t$,

$$\tilde{\mathbf{P}}(M_{N_t} - M_0 < -\epsilon t) \leq \tilde{\mathbf{P}}(N_t > \beta t) + \tilde{\mathbf{P}}\left(\min_{j \leq \beta t} M_j - M_0 < -\epsilon t\right)$$

$$\leq e^{-Aht} + e^{-b(\epsilon^2 t^2)/(\beta t)} \tag{20}$$

$$\leq e^{-d\epsilon^2 t}, \tag{21}$$

where $d > 0$ is a universal constant. Note that $M_0 = Y_0$ and $|M_{N_t} - Y_{t \wedge W}| \leq \frac{1}{2}$ for all $t$. It follows that for every $\epsilon > 0$ we have

$$\tilde{\mathbf{P}}(Y_{t \wedge W} - Y_0 < -\epsilon t) \leq e^{-d\epsilon^2 t} \tag{22}$$

for sufficiently large $t$, if we incorporate an extra factor of, say, $\frac{1}{2}$ into $d$.

To complete the proof of Lemma 10, note that

$$\tilde{\mathbf{P}}(W > t) \leq \tilde{\mathbf{P}}\left(B_t < \frac{Ct}{8(\rho+1)}\right) + \tilde{\mathbf{P}}\left(W > t, B_t \geq \frac{Ct}{8(\rho+1)}\right)$$

$$\leq e^{-ct/(\rho+1)^2} + \tilde{\mathbf{P}}\left(W > t, B_t \geq \frac{Ct}{8(\rho+1)}\right), \tag{23}$$

where the second inequality uses equation (18). But if $W > t$, then $Y_t$ is $Y_{t \wedge W}$, so

$$\tilde{\mathbf{P}}\left(W > t, B_t \geq \frac{Ct}{8(\rho+1)}\right)$$

$$\leq \tilde{\mathbf{P}}\left(Y_{t \wedge W} \leq U_t + \frac{4t\alpha}{\tilde{\rho}} - \frac{\alpha Ct}{8(\rho+1)}\right)$$

$$= \tilde{\mathbf{P}}\left(Y_{t \wedge W} \leq U_t - \frac{\alpha Ct}{16(\rho+1)}\right)$$

$$\leq \tilde{\mathbf{P}}\left(Y_{t\wedge W} - Y_0 \leq -\frac{\alpha C t}{16(\rho+1)} - Y_0\right),$$

where the equality holds because $\tilde{\rho} = 64C^{-1}(\rho+1)$ and the second inequality holds because $U_t \leq 0$. But for sufficiently large $t$, we have $|Y_0| \leq \frac{\alpha C t}{32(\rho+1)}$, so

$$\tilde{\mathbf{P}}\left(W > t, B_t \geq \frac{Ct}{8(\rho+1)}\right) \leq \tilde{\mathbf{P}}\left(Y_{t\wedge W} - Y_0 \leq -\frac{\alpha C t}{32(\rho+1)}\right)$$

$$\leq \exp\left\{\frac{-d\alpha^2 C^2 t}{(32)^2(\rho+1)^2}\right\},$$

where the second inequality uses (22). Combining this with (23) completes the proof. □

## APPENDIX

The following lemma was used in the proof of Lemma 10:

LEMMA A.1. *Let $X_t$ be an irreducible continuous-time Markov chain on a finite state space $S$. Fix $s \geq 0$, and for $x \in S$, let $W_x = |\{0 < t < s : X_t = x\}|$ be the amount of time the chain spends in state $x$ up to time $t$. Then there is a constant $c$, independent of $s$, such that for any event $A \in \sigma(W_x : x \in S)$, and any $y \in S$, we have*

$$\mathbf{P}_y(A) \leq c \max_{z \in S} \tilde{\mathbf{P}}_z(A),$$

*where $\tilde{\mathbf{P}}(A)$ denotes the probability of $A$ when we run the time reversal of $X_t$.*

PROOF. Suppose that $X_t$ has stationary distribution $\pi$ and let $c = |S| \times \max_{z,y} \pi(z)/\pi(y)$. Suppose that $X_0$ is distributed according to $\pi$. Then $\{X_t : 0 \leq t \leq s\}$ has the same distribution as $\{\tilde{X}_{s-t} : 0 \leq t \leq s\}$, where $\tilde{X}$ is the time reversal of $X$, with $\tilde{X}_0$ having distribution $\pi$.

Fix an event $A \in \sigma(W_x : x \in S)$. For states $i, j$, let $A_{i,j} = A \cap [X_0 = i, X_s = j]$. Since $\{X_t : 0 \leq t \leq s\}$ has the same distribution as $\{\tilde{X}_{s-t} : 0 \leq t \leq s\}$, it is clear that $\mathbf{P}(A_{i,j}) = \tilde{\mathbf{P}}(A_{ji})$, for all $i$ and $j$. Hence,

$$\mathbf{P}_y(A) = \frac{1}{\pi(y)} \mathbf{P}\left(\bigcup_z A_{y,z}\right)$$

$$\leq \frac{|S|}{\pi(y)} \max_z \mathbf{P}(A_{y,z})$$

$$= \frac{|S|}{\pi(y)} \max_z \tilde{\mathbf{P}}(A_{z,y})$$

$$\leq \frac{|S|}{\pi(y)} \max_z \pi(z) \tilde{\mathbf{P}}(A)$$

$$\leq c \max_z \tilde{\mathbf{P}}_z(A). \qquad \square$$

The following proposition was used in the proof of Lemma 6:

PROPOSITION A.2. *Let $X_t$ be the continuous-time simple symmetric random walk on the integers in which moves in each direction are made at rate 1. Then there is a universal constant $b > 0$ such that for every $r \geq 1$, we have $\mathbf{P}_0(X_t \neq 0$ for all $t \in [1, r^2]) \geq br^{-1}$.*

PROOF. Let $R = \lceil r \rceil$. Since the number of jumps by time 1 is Poisson(2), we have

$$\mathbf{P}_0(X_1 = 1) \geq \tfrac{1}{2}\mathbf{P}(\text{exactly one jump by time } 1) = e^{-2},$$

and given that $X_1 = 1$, the conditional probability of hitting $R$ before returning to 0 is $\frac{1}{R}$ since $X_t$ is a martingale. Call this event $B$.

It is well known that for any $t \geq 0$, the probability that simple random walk deviates to the left of its starting point by more that $\sqrt{t}$ up to time $t$ can be bounded away from 0 (uniformly in $t$). Hence, given $B$, the conditional probability that $X_t$ does not return to 0 by time $R^2$ is bounded away from 0. It follows that

$$\mathbf{P}_0(X_t \neq 0 \text{ for all } t \in [1, R^2]) \geq bR^{-1}$$

for a universal constant $b$. Since $R \in [r, 2r]$, the result follows if we incorporate an extra factor of $\frac{1}{2}$ into $b$. $\square$

**Acknowledgments.** I would like to thank P. Diaconis, S. Evans, J. Gravner, R. Lyons, F. Martinelli, Y. Peres and J. Quastel for useful discussions.

Fabio Martinelli (following up on a suggestion by Oded Schramm) showed me how to prove the optimality of my upper bounds, as did one of the referees. I want to thank Fabio Martinelli and one of the referees for bringing this to my attention.

## REFERENCES

[1] CAPUTO, P. (2004). Spectral gap inequalities in product spaces with conservation laws. In *Stochastic Analysis on Large Interacting Systems* (T. Funaki and H. Osada, eds.) 53–88. Math. Soc. Japan, Tokyo. MR2073330
[2] DIACONIS, P. and SALOFF-COSTE, L. (1993). Comparison theorems for reversible Markov chains. *Ann. Appl. Probab.* **3** 696–730. MR1233621
[3] DIACONIS, P. and STROOCK, D. (1991). Geometric bounds for eigenvalues of Markov chains. *Ann. Appl. Probab.* **1** 36–61. MR1097463

[4] HORN, R. A. and JOHNSON, C. R. (1985). *Matrix Analysis.* Cambridge Univ. Press. MR0832183

[5] LANDIM, C., SETHURAMAN, S. and VARADHAN, S. (1996). Spectral gap for zero-range dynamics. *Ann. Probab.* **24** 1871–1902. MR1415232

[6] LIGGETT, T. M. (1989). Exponential $L_2$ convergence of attractive reversible nearest particle systems. *Ann. Probab.* **17** 403–432. MR0985371

[7] LU, S. and YAU, H.-T. (1993). Spectral gap and logarithmic Sobolev inequality for Kawasaki and Glauber dynamics. *Comm. Math. Phys.* **156** 399–433. MR1233852

[8] MCDIARMID, C. (1989). On the method of bounded differences. In *Surveys in Combinatorics* (J. Siemons, ed.) 148–188. Cambridge Univ. Press. MR1036755

[9] SINCLAIR, A. J. (1992). Improved bounds for mixing rates of Markov chains and multicommodity flow. *Combin. Probab. Comput.* **1** 351–370. MR1211324

[10] SALOFF-COSTE, L. (1997). Lectures on finite Markov chains. *Lecture Notes in Math.* **1665** 301–413. Springer, Berlin. MR1490046

[11] VARADHAN, S. R. S. (1994). Nonlinear diffusion limit for a system with nearest neighbor interactions. II. In *Asymptotic Problems in Probability Theory*: *Stochastic Models and Diffusion on Fractals* (K. D. Elworthy and N. Ikeda, eds.) 75–128. Longman Sci. Tech., Harlow. MR1354152

[12] WILSON, D. (2004). Mixing times of lozenge tiling and card shuffling Markov chains. *Ann. Appl. Probab.* **14** 274–325. MR2023023

DEPARTMENT OF MATHEMATICS
UNIVERSITY OF CALIFORNIA
DAVIS, CALIFORNIA 95616
USA
E-MAIL: morris@math.ucdavis.edu